\tolerance=10000
\raggedbottom

\baselineskip=15pt
\parskip=1\jot

\def\sk{\vskip 3\jot}

\def\heading#1{\vskip3\jot{\noindent\bf #1}}
\def\label#1{{\noindent\it #1}}


\def\ref#1;#2;#3;#4;#5.{\item{[#1]} #2,#3,{\it #4},#5.}
\def\refinbook#1;#2;#3;#4;#5;#6.{\item{[#1]} #2, #3, #4, {\it #5},#6.} 
\def\refbook#1;#2;#3;#4.{\item{[#1]} #2,{\it #3},#4.}


\def\({\bigl(}
\def\){\bigr)}

\def\Ex{{\rm Ex}}
\def\Var{{\rm Var}}

\def\la{\lambda}
\def\ta{\tau}

\def\sumin{\sum_{1\le i\le n}}
\def\avein{{1\over n}\sumin}
\def\bpi{b_{\pi(i)}}
\def\bti{b_{\ta(i)}}

{
\pageno=0
\nopagenumbers
\rightline{\tt egp.i.arxiv.tex}
\vskip1in

\centerline{\bf A Bound on the Variance of the Waiting Time in a Queueing System}
\vskip0.5in

\centerline{Patrick Eschenfeldt}
\centerline{\tt peschenfeldt@hmc.edu}
\sk

\centerline{Ben Gross}
\centerline{\tt bgross@hmc.edu}
\sk

\centerline{Nicholas Pippenger}
\centerline{\tt njp@math.hmc.edu}
\sk

\centerline{Department of Mathematics}
\centerline{Harvey Mudd College}
\centerline{1250 Dartmouth Avenue}
\centerline{Claremont, CA 91711}
\vskip0.5in

\noindent{\bf Abstract:}
Kingman has shown, under very weak conditions on the interarrival- and sevice-time distributions, that First-Come-First-Served minimizes the variance of the waiting time
among possible service disciplines.
We show, under the same conditions, that Last-Come-First-Served maximizes the variance of the waiting time, thereby giving an upper bound on the variance among all disciplines.
\vskip0.5in
\leftline{{\bf Keywords:} Queueing theory, waiting time, service discipline.}
\sk
\leftline{{\bf Subject Classification:} 60K26, 90B22}
\vfill\eject
}

\heading{1.  Introduction}

It is well known that the average waiting time $\Ex[W]$ in a queueing system does not depend on the service discipline (that is, on the rule specifying which waiting customer is to be served when a server becomes free).
This fact may be seen from  the classic result of Little [L] that the average number $\Ex[K]$ of customers in the system is equal to the average rate $\la$ at which customers arrive multiplied by the average sojourn time $\Ex[S]$
(where the sojourn time $S$ is the sum of the waiting time $W$  and service time $V$); 
since the service discipline does not affect $\Ex[K]$, $\la$ or $\Ex[V]$, it cannot affect $\Ex[W]$.
(Little assumes that the various processes involved are stationary and ergodic, and that the expectations
in question are finite.)
The variance $\Var[W]$ of the waiting time, however, {\it does\/} depend on the service discipline,
and Kingman [K] has shown that  ``first-come-first-served'' (FCFS) minimizes this variance.
(Kingman assumes that the ``null state'' of the empty queue is recurrent.
His result does not require that the variance of the waiting time be finite, if one interprets
it as saying that if the variance for FCFS is infinite, then so is that for any other service discipline.)
Our goal in this paper is to show that ``last-come-first-served'' (LCFS) maximizes the variance.
This result provides an upper bound on the variance of any service discipline.
(Our result is an elaboration of  Kingman's, and holds under the same assumptions.
It shows that if the variance for any service discipline is infinite, then so is that for LCFS.
For the case of the $M/M/s$ system, Vaulot [V] has derived the waiting time distribution for
LCFS. 
For $M/M/1$, Riordan [R] gives explicit expressions for the second moments (conditional on $W>0$) for both FCFS, where $\Ex[W^2\mid W>0] =  2/(1-\la)^2$, and LCFS, where
$\Ex[W^2\mid W>0] =  2/(1-\la)^3$.
Combined with the results $\Pr[W>0] = \la$ and $\Ex[W] = \la/(1-\la)$,
which are independent of the service discipline,
these expressions yield $\Var[W] = \la(2-\la)/(1-\la)^2$ for FCFS and
$\Var[W] = \la(2-\la+\la^2)/(1-\la)^3$ for LCFS.)

Kingman's argument proceeds as follow.
Since $\Var[W] = \Ex[W^2] - \Ex[W]^2$ and $\Ex[W]$ is not affected by the service discipline,
the discipline that minimizes $\Ex[W^2]$ also minimizes $\Var[W]$.
Let $a_1 < a_2 < \cdots < a_n$ be the arrival times of customers during a busy period,
and let $b_1 < b_2 < \cdots < b_n$ be the beginnings of the service intervals during this busy period.
Suppose that, for $1\le i\le n$, the customer arriving at time $a_i$ is served at time $b_{\pi(i)}$,
where $\pi$ is a permutation of $\{1, \ldots, n\}$ that depends on the service discipline.
Since the waiting time of customer $i$ is $\bpi - a_i$,
the contribution to $\Ex[W^2]$ from this busy period is 
$$\eqalign{
\avein (\bpi - a_i)^2 &= \avein \bpi^2 - {2\over n}\sumin a_i\bpi + \avein a_i^2 \cr
&= \avein b_i^2 - {2\over n}\sumin a_i\bpi + \avein a_i^2. \cr
}$$
Thus the discipline that maximizes $A(\pi) = \sumin a_i\bpi$ minimizes $\Ex[W^2]$.
The ``rearrangement inequality'' (see, for example, Hardy, Littlewood and P\'{o}lya [H],
Theorem 368)
says that $A(\pi)$ is maximized by the identity permutation
$\pi(i)=i$, which corresponds to FCFS.

This argument shows that to maximize $\Var[W]$, we should minimize $A(\pi)$.
The rearrangement inequality says that, among all permutations $\pi$, $A(\pi)$ is minimized
by the reverse permutation, $\pi(i) = n+1-i$.
But LCFS does not always give this permutation!
If $i<j$, then customer $i$ will be served before customer $j$ in FCFS, but customer $i$ is not be necessarily after customer $j$ in LCFS (because customer $j$ might not even have arrived when customer $i$ is ready to be served). 
Indeed, customer $1$, whose arrival initiates the busy period, is served immediately by all service
disciplines.
Thus to maximize $\Var[W]$, we must minimize $A(\pi)$, not among all permutations $\pi$,
but rather among those in the set $\Pi$ of permutations that can be realized by a service discipline.
This set  depends on the interleaving of the $a_1 < a_2 < \cdots < a_n$ among the 
$b_1 < b_2 < \cdots < b_n$ in the particular busy period under study.
In the next section we shall study this situation, and show that $A(\pi)$ is indeed
minimized over $\pi\in\Pi$ by the permutation that is realized  by LCFS.
\sk

\heading{2. The Minimizing Property of LCFS}

We have observed that for any $\pi\in\Pi$, we have $a_1 = b_1$ and $\pi(1)=1$ so the term 
$a_1 b_{\pi(1)} = a_1 b_1$ in $A(\pi)$ is independent of $\pi$.
Thus we may restrict our attention to the restriction of $\pi$ to the set $\{2, \ldots, n\}$.
We shall assume that no two of the numbers $a_2, \ldots, a_n, b_2,\ldots, b_n$ are equal.
(This event occurs with probability one if the arrival times and service times have absolutely continuous distributions, and a trite perturbation argument shows that it entails no loss of generality
in other cases.)
Then
$$\Pi = \big\{\pi\hbox{\ a\ permutation\ of\ }\{1, \ldots, n\}: 
\pi(1)=1 \hbox{\ and\ }a_i < \bpi \hbox{\ for\ } 2\le i\le n\big\}.$$

Let $\ta$ be the permutation realized by LCFS, and let $\pi\in\Pi$. 
We shall show that
$$A(\ta)\le A(\pi). \eqno(2.1)$$
Say that a pair $(i,j)$ is a {\it bad\/} for $\pi$ if $a_i < a_j < b_{\pi(i)} < b_{\pi(j)}$.
The pair $(i,j)$ being bad for $\pi$ signifies that the discipline realizing $\pi$, confronted with customers $i$ and $j$ who could be served in either order, served first the one who arrived earlier.
Thus $\ta$ is the unique permutation in $\Pi$ that has no bad pairs.
This observation proves (2.1) for permutations $\pi\in\Pi$ that have no bad pairs.

It remains to prove (2.1) for permutations $\pi\in\Pi$ that  have one or more bad pairs.
Suppose, to obtain a contradiction, that $\pi$ is a counterexample; that is (1) $\pi\in\Pi$, (2)
$\pi$ has one or more bad pairs, and (3)
$$A(\pi) < A(\ta). \eqno(2.2)$$
We may assume that $\pi$ is a ``smallest possible''  counterexample to (2.1); that is, that (a)
the busy period under study has the smallest possible value of $n$, and (b) among counterexamples for that value of $n$,
$\pi$ has the smallest possible number of bad pairs.

Let $n-1$ left parentheses ``$($'' be positioned at the points $a_2, \ldots, a_n$ on the real line,
and let $n-1$ right parentheses ``$)$'' be positioned at the points $b_2, \ldots, b_n$.
Since $\bpi > a_i$ for all $2\le i\le n$, these $n-1$ pair of parentheses will match in the usual way, and the left parenthesis at $a_i$ will be matched with the right parenthesis at $\bti$ (right parentheses match left parentheses on a LCFS basis).

Since the sequence of parentheses begins with a left parenthesis and ends with a right parenthesis,
there must at some point be a left parenthesis (say at $a_k$) that immediately followed by a right parenthesis (say at $b_l$), and we must have $\ta(k)=l$.
If we also have $\pi(k)=l$, then we may remove $a_k$ and $b_l$ from their respective lists; this will not remove any bad pairs, so we thereby obtain a counterexample with a smaller value of $n$.
Thus we must have $\pi(k)\not=l$.
Let $i$ be such that $\pi(i)=l$.
We must then have $i<k$ and $\pi(i)<\pi(k)$, so that the pair $(i,k)$ is bad for $\pi$.
Define a new permutation $\pi'\in\Pi$ by swapping the values assigned by $\pi$  to $i$ and $k$; 
that is
$$\pi'(m) = \cases{
\pi(i), &if $m=k$, \cr
\pi(k), &if $m=i$, \cr
\pi(m), &otherwise. \cr
}$$
The pair $(i,k)$, which is bad for $\pi$,  is not bad for $\pi'$.
If there are any $j$ in the range $i < j < k$ for which $\pi(j)$ satisfies $\pi(i) < \pi(j) < \pi(k)$, then
the pairs $(i,j)$ and $(j,k)$, which are bad for $\pi$, are not bad for $\pi'$.
But no pairs that are not bad for $\pi$ are bad for $\pi'$.
Thus $\pi'$ has fewer bad pairs than $\pi$, so we must have
$$A(\ta)\le A(\pi'), \eqno(2.3)$$
else $\pi'$ would be a counterexample with fewer bad pairs than $\pi$.
But the inequalities $a_i < a_k < \bpi < b_{\pi(k)}$ imply $a_i b_{\pi'(i)} + a_k b_{\pi'(k)} < 
a_i b_{\pi(i)} + a_k b_{\pi(k)}$, which in turn implies
$$A(\pi') < A(\pi) \eqno(2.4)$$
(because these are the only two terms that differ between $A(\pi')$ and $A(\pi)$).
Combining (2.3) and (2.4) contradicts (2.2), and thus completes the proof of (2.1).
\sk

\heading{3. Conclusion}

We have shown that, among all serve disciplines, LCFS maximizes the variance of the waiting time.
It is not hard to see that a simple modification of our proof (reversing some of the inequalities)
yields Kingman's result that FCFS minimizes the variance of the waiting time.
Finally, another modification (observing that the swap $\pi\mapsto\pi'$ preserves the sum of the waiting times)
yields the result that the average waiting time is not affected by the service discipline.
\sk

\heading{4. Acknowledgment}

The research reported here was supported
by Grant CCF  0917026 from the National Science Foundation.
\sk

\heading{5. References}

\refbook H; G. H. Hardy, J. E. Littlewood and G. P\'{o}lya;
Inequalities;
Cambridge University Press, London, 1934.

\ref K; J. F. C. Kingman;
``The Effect of Queue Discipline on Waiting Time Variance'';
Math.\ Proc.\ Cambridge Phil.\ Soc.; 58:1 (1962) 163--164.

\ref L; J. D. C. Little;
``A Proof for the Queuing Formula: $L = \la W$'';
Oper.\ Res.; 9 (1961) 383--387.

\refbook R; J. Riordan;
Stochastic Service Systems;
John Wiley and Sons, New York, 1962.


\ref V; \'{E}. Vaulot;
``D\'{e}lais d'attente des appels t\'{e}l\'{e}phoniques dans l'ordre inverse de leurs arriv\'{e}e'';
Comptes Rendus Acad.\ Sci.\ Paris; 238 (1954) 1188--1189.

\bye